\documentclass[11pt]{amsart}
\usepackage{geometry}                % See geometry.pdf to learn the layout options. There are lots.
\usepackage{graphicx}
\usepackage{amssymb}
\usepackage{epstopdf}
\usepackage{parskip}

\title{ON THE RIBENBOIM - LANDAU CONJECTURE}
\author{F. Sidokhine}
\date{\today}                                           % Activate to display a given date or no date
\begin{document}

\begin{abstract}
We prove the Ribenboim hypothesis, which states that if, starting from some integer \(N\), consecutive prime numbers \(p_ {n}\), \(p_{n+1}\) satisfy the inequality \(\sqrt {p_ {n+1}}-\sqrt{p_{n}} <1\), then the Landau problem \# 4 (1912) has a positive solution.
\end{abstract}

\maketitle
\section{Introduction}

P. Ribenboim,  [1],  p.191, (2004), hypothesizes the solution of the Landau problem using the inequality \(\sqrt{p_{n+1}}-\sqrt{p_{n}}<1\), formulating this problem as follows: "\textit{from this inequality, if true, it would follow that between the squares of any two consecutive integers, there is always a prime}" (the Landau Problem \# 4 (1912), [2]). In this note, we will prove the Ribenboim conjecturel, and also consider some other problems with prime numbers.

\section{The Landau - Ingham problems and the Ribenboim conjecture }

Before proving the Ribenboim conjecture, we will discuss the Ribenboim hypothesis, in the context of understanding what it states in terms of the Landau problem itself. In other words, let the Landau hypothesis be true, then:

Theorem 1. Let there be a positive integer constant \textit{N}, such that for any integer \textit{n} greater than \textit{N} between squares of consecutive integers \textit{}\(n\) and \textit{}\(n+1\), there is always a prime, then the consecutive prime numbers \(p_{m}, p_{m+1}\) satisfy the relation \(\sqrt{p_{m+1}}-\sqrt{p_{m}} = O(1\)).
\begin{proof}
There is a positive integer\textit{ M} such that \( p_{M-1}<\left(N+1\right)^{2}<p_{M}\), where \( p_{M-1}\), \( p_{M}\) are primes. Let \( C\) be an arbitrary integer greater than 2, and there exists an integer \( m_{0}\) greater than \textit{M}, such that primes \( p_{m_{0}}, p_{m_{0}+1}\) satisfies the inequality \( \sqrt{p_{m_{0}+1}}-\sqrt{p_{m_{0}}}>C\). We have two possibilities for \(\sqrt{p_{m_{0}+1}}\), or \(\left[\sqrt{p_{m_{0}}}\right] <\sqrt{p_{m_{0}}} <\sqrt{p_{m_{0}+1}} <\left[\sqrt{p_{m_{0}}}\right]+1\), or \(\left[\sqrt{p_{m_{0}}}\right]+1<\sqrt{p_{m_{0}+1}} <\left[\sqrt{p_{m_{0}}}\right]+2\), since \(\left[\sqrt{p_{m_{0}}}\right]>N\). Thus, \(\sqrt{p_{m_{0}+1}}-\sqrt{p_{m_{0}}}<2\){\small  } and we have a contradiction that proves Theorem 1.
\end{proof}

Thus, the Landau hypothesis states that, the consecutive prime numbers \( p_{m},  p_{m+1}\), where \( m\) is greater than \( M\), must satisfy the inequality \(\sqrt{p_{m+1}}-\sqrt{p_{m}}<C\), where \( C\) is some positive constant. The question arises at what value of the constant \( C\) the converse statement is true. The Ribenboim hypothesis states that this value of the constant \( C\) is 1.

Theorem 2 (Ribenboim's conjecture). Let there be a positive integer constant \textit{M}, such that for any integer \textit{m} greater than \textit{M} the consecutive prime numbers \(p_{m}, p_ {m+1}\) satisfy the inequality \(\sqrt{p_ {m+1}}-\sqrt{p_{m}}<1\), then between squares of consecutive integers \(n\) and \textit{}\(n+1\), where \( n>p_{M}\), there is always a prime number.

\begin{proof}

 To prove Theorem 2, the following Lemmas 1 and 2 are necessary:

Lemma 1. There is a positive integer constant \textit{M} greater than 2, such that for any integer \( m\) greater than \textit{M} the open interval of real numbers \(\left(p_{m+1}-2\sqrt{p_{m+1} },p_{m+1}\right)\), where \( p_{m+1}\) is a prime, contains a prime number.
\begin{proof}
According to the condition of Theorem 2, there exists an integer \textit{M} such that for any integer \(m\) greater than \textit{M}, the inequality \(p_{m+1}-p_{m}<\sqrt{p_ {m+1}}+\sqrt{p_{m}}<2\sqrt{p_{m+1}}\) holds, where \(p_ {m}, p_{m+1}\) are consecutive prime numbers. So we have \(p_{m+1}-2\sqrt{p_{m+ 1}}<p_{m}<p_{m+1}\), and Lemma 1 is proved.
\end{proof}

Lemma 2. Let there be a positive integer constant \textit{M} greater than 2, such that for any  integer \textit{m }greater than \textit{M} the difference of consecutive prime numbers \( p_{m}, p_{m+1}\) satisfies the inequality \( p_{m+1}-p_{m}<2\sqrt{p_{m+1}}\), then for any real number \( x\) greater than \( p_{M+1}\), the open interval of real numbers \( (x-2\sqrt{x}, x)\) contains a prime number.
\begin{proof}
 Let there be some real number \( x_{0}>p_{M+1}\) such that the open interval \( (x_{0}-2\sqrt{x_{0}}, x_{0})\)  does not contain primes. Let us construct an open interval \(\left(p_{k+1}-2\sqrt{p_{k+1} },p_{k+1}\right)\), where \( p_{k}<x_{0}<p_{k+1}\), since \( k>M\), therefore, the open interval \(\left(p_{k+1}-2\sqrt{p_{k+1} },p_{k+1}\right)\) contains a prime, according to Lemma 1. However, on the other hand, the open interval \(\left(p_{k+1}-2\sqrt{p_{k+1} },p_{k+1}\right) =\left(p_{k+1}-2\sqrt{p_{k+1} },x_{0}\right)\cup\left[x_{0},p_{k+1}\right)\), where \( p_{k+1}-2\sqrt{p_{k+1} }<x_{0}\) otherwise \( x_{0}<p_{k}<p_{k+1}\), \(\left(p_{k+1}-2\sqrt{p_{k+1} },x_{0}\right)\subset\left(x_{0}-2\sqrt{x_{0}}, x_{0}\right)\),\textit{ }since \( x_{0}-2\sqrt{x_{0}}<p_{k+1}-2\sqrt{p_{k+1} }\), does not contain prime numbers. We have obtained a contradiction, and Lemma 2 is proved. 
\end{proof}

Let the integer \( m_{0}\) greater than \( p_{M}\) and such that the interval \( (m_{0}^{2}, (m_{0}+1)^{2})\) does not contain primes. Using directly Lemma  2, where \( x = (m_{0}+1)^{2}\), we conclude that the interval \( (m_{0}^{2}-1, (m_{0}+1)^{2})\) contains a prime number. Thus, we have \( m_{0}^{2}-1<m_{0}^{2}<p<(m_{0}+1)^{2}\), where \textit{p} is the prime number, so as \( m_{0}^{2}-1\) and \( m_{0}^{2}\) are neighboring composite integers. We have obtained a contradiction, and Theorem 2 is proved.
\end{proof}

Note. Numerical tests (see [1], p. 191) confirm that, starting from some fixed integer \(N\), consecutive primes \(p_ {n}\), \(p_{n+1}\) satisfy the inequality \(\sqrt{p_{n+1}}-\sqrt{p_{n}} <1\), and there is every reason to believe that this hypothesis (Andrica's hypothesis) is true with high probability.

Developing the Ribenboim approach, we linked the solution of some problems with prime numbers with the relation \(\sqrt[n]{p_{m+1}} -\sqrt[n]{p_{m}} = O(1)\){\small}, where \(n\) and \(m\) are integers and  \(n>\) 1. Here are some examples:

Theorem 3 (direct). Let there be a  positive  integer constant  \textit{N}, such that for any integer \textit{n} greater than \textit{N}  between  squares of any two rational numbers \( n\), \( n+\frac{1}{2}\) and \(  n+\frac{1}{2}\),\textit{ }\( n+1\), there are always prime numbers, then there exists a positive  integer \textit{M}, such that for any integer\textit{ m }greater than \textit{M} the consecutive prime numbers \( p_{m},  p_{m+1}\) satisfy the relation \(\sqrt{p_{m+1}}-\sqrt{p_{m}}<1\).

Theorem 4 (inverse). Let there be a positive integer constant \textit{M}, such that for any integer \textit{m }greater than \textit{M}  consecutive primes \( p_{m},  p_{m+1}\) satisfy the inequality \(\sqrt{p_{m+1}}-\sqrt{p_{m}}<\frac{1}{2}\), then between squares of any two rational  numbers \( n\), \( n+\frac{1}{2}\) and \(  n+\frac{1}{2}\),\textit{ }\( n+1\), where the integer \( n>p_{M}\), there are always prime numbers.

Theorem 5 (direct). Let there be a positive  integer \textit{N}, such that for any integer \textit{n} greater than \textit{N} between  cubes of  consecutive integers\textit{ }\( n\) and\textit{ }\( n+1\), there is always a prime, then consecutive prime numbers \( p_{m},  p_{m+1}\) satisfies the relation \(\sqrt[3]{p_{m+1}}   -\sqrt[3]{p_{m}} = O(1)\).

Theorem 6 (inverse, Ingham's problem). Let there be a positive integer \textit{M}, such that for any integer \textit{m} greater than \textit{M} the primes \( p_{m},  p_{m+1}\) satisfy the inequality \(\sqrt[3]{p_{m+1}} - \sqrt[3]{p_{m}} <1-\delta\), where \( \delta\) is any real number satisfying the inequality \( 0<\delta <1\), then between cubes of consecutive integers\textit{ }\( n\) and\textit{ }\( n+1\), where \( n>p_{M}\), there is always a prime number.
\begin{proof}
 Ingham's problem implies that consecutive prime numbers \(p_{n}\), \(p_{n+1}\) satisfy the relation \(\sqrt{p_{n+1}}-\sqrt{p_{n}}=O(p_{n+1}^{1/6})\). According to [3], the relation \(\sqrt{p_{n+1}}-\sqrt{p_{n}} = O(p_{n+1}^{3/46})\) is true, hence \(\sqrt{p_{n+1}}-\sqrt{p_{n}}=o(p_{n+1}^{1/6})\) and consecutive primes \(p_{n}\), \(p_{n+1}\) satisfy the relations \({p_{n+1}}-{p_{n}} = o(p_{n+1}^{2/3})\) and \(p_{n+1}^{1/3}\) \( -\) \( p_{n}^{1/3}\) \(=o(1)\). Thus, starting from some integer \(N\), consecutive prime numbers \(p_{n}\), \(p_{n+1}\)  must satisfy the inequality \(p_{n+1}^{1/3}\) \( -\) \( p_{n}^{1/3}\) \(<1-\delta\).

To prove Theorem 6, the following Lemmas 3, 4 and 5 are necessary:

Lemma 3. Let there be a positive integer constant \textit{M}, such that for any integer \( m\) greater than \textit{M} consecutive prime numbers \( p_{m}\), \( p_{m+1}\) satisfy the inequality \(p_{n+1}^{1/3}\) - \(p_{n}^{1/3}\) \(<1-\delta\), then consecutive primes \( p_{m}\), \(p_{m+1}\) satisfy the inequality \(  p_{m+1}-p_{m}<3(1-\delta )p_{m+1}^{2/3}\).

Lemma 4. Let there be a positive integer \textit{M}, such that for any integer \( m\) greater than \textit{M} the open interval {\small (}\( p_{m+1}-3(1-\delta )p_{m+1}^{2/3},p_{m+1})\), where \( p_{m+1}\) is a prime and \(\sqrt[3]{p_{M}}>3(1-\delta )\), contains a prime number.

Lemma 5. Let there be a positive integer constant \textit{M}, such that for any integer \textit{m }greater than \textit{M} consecutive primes \( p_{m},  p_{m+1}\) satisfy the inequality \(  p_{m+1}-p_{m}<3(1-\delta )p_{m+1}^{2/3}\), where \(\sqrt[3]{p_{M}}>3(1-\delta )\), then for any real number \( x\) greater than \( p_{M}\), the open interval of real numbers \( (x-3(1-\delta )x^{2/3}, x)\) contains a prime number.

Let the integer \(m_{0}\) be greater than \(p_{M}\) and such that the interval \( (m_{0}^{3} ,(m_{0}+1)^{3})\) does not contain prime numbers. Using Lemma 5, where \(x =\)\((m_{0}+1)^{3}\), we conclude that the interval \(((m_{0}+1)^{3}-3(1-\delta )(m_{0}+1)^{2}, (m_{0}+1)^{3})\) contains a prime number. Then we have

\begin{center}
\( m_{0}^{3}<(m_{0}+1)^{3}-3(1-\delta )(m_{0}+1)^{2}\) \( <p<(m_{0}+1)^{3}\).
\end{center}

This inequality is true since \( m_{0}>p_{M}\), where \( p_{M}\) is also greater than  max\((27\left(1-\delta\right)^{3}, 3/\delta )\). Thus, we have obtained a contradiction, and Theorem 6 is proved.
\end{proof}

\section{Solving Problems with Primes Numbers within the framework of the Ribenboim hypothesis}

Using Ribenboim's hypothesis we can prove stronger results, Theorems 7, 9 and 10.

Theorem 7. Let the consecutive primes \( p_{n-1},  p_{n}\) satisfy the relation \(\sqrt{p_{n}}-\sqrt{p_{n-1}} = O(1)\), then there exists a positive integer constant \( M\), such that between  cubes of consecutive integers \(m\)\textit{ }and\textit{ }\(  m+1\),\textit{ }where  \( m > M\), there is always a prime number.

Theorem 7 is a direct consequence of the following statement:

Theorem 8. Let there be positive integer constants \( A\) and \textit{B}, such that for any integer \( n\) greater than \( A\) the consecutive primes \( p_{n-1},  p_{n}\) satisfy the inequality \(\sqrt{p_{n}}-\sqrt{p_{n-1}}<B/2,\) where \(  B\geq 4\); \(\sqrt{p_{A}}>B\), then between  cubes of any two real numbers \( x-1\)\textit{ }and\textit{ }\(  x\),\textit{ }where \textit{x}\( \) greater than \( p_{A}\), there is always a prime number.
\begin{proof}

To prove Theorem 8, the following Lemmas 6, 7 and 8 are necessary:

Lemma 6. There be positive integer constants \textit{A} and \textit{B}, such that for any integer \( n\) greater than \textit{A} consecutive prime numbers \( p_{n-1}, p_{n}\) satisfy the inequality \(\sqrt{p_{n}}-\sqrt{p_{n-1}}<B/2\), where the integer \( B\geq 4\), then the difference of consecutive prime numbers \( p_{n-1},  p_{n}\) satisfies the inequality \( p_{n}-p_{n-1}<B\sqrt{p_{n}}\).

Lemma 7. Let there be positive  integers \textit{A} and \textit{B}, such that for any integer \( n\) greater than \textit{A} the interval \(\left(p_{n}-B\sqrt{p_{n} },p_{n}\right)\), where \( p_{n}\) is a prime and \(\sqrt{p_{A}}>B\), contains a prime.

Lemma 8. Let there be positive  integers \textit{A} and \textit{B}, such that for any  integer \textit{n }greater than \textit{A} the difference of consecutive primes \( p_{n-1},  p_{n}\) satisfies the inequality \( p_{n}-p_{n-1}<B\sqrt{p_{n}}\), where \(\sqrt{p_{A}}>B\), then for any real number \( x\) greater than \( p_{A}\), the open interval of real numbers \( (x-B\sqrt{x}, x)\) contains a prime number.

Let the real number \( x_{0}\) greater than \( p_{A}\), such that the interval \( ((x_{0}-1)^{3},x_{0}^{3})\) does not contain prime numbers. Using Lemma 8, where \( x =\) \( x_{0}^{3}\), we conclude that the interval \( (x_{0}^{3}-Bx_{0}^{3/2}, x_{0}^{3})\) contains a prime number. Then we have

\begin{center}
\( (x_{0}-1)^{3}<x_{0}^{3}\) – 3 \( x_{0}^{2}+3x_{0}<x_{0}^{3}-Bx_{0}^{3/2}\) as a consequence of the inequality  \(\frac{B}{3\sqrt{x_{0}}} +\frac{1}{x_{0}}<1\).
\end{center}

This inequality is true since \( x_{0}>p_{A}>B^{2}>4\), so \( (x_{0}-1)^{3}< x_{0}^{3}-Bx_{0}^{3/2}<p<x_{0}^{3}\), where \( p\) is some prime. Thus, we have obtained a contradiction, and Theorem 8 is proved.
\end{proof}

Theorem 9. Let there be positive  integer constants \textit{A} and \textit{B}, such that for any  integer \( n\) greater than \( A\) consecutive primes \( p_{n-1},  p_{n}\) satisfy the inequality \(\sqrt{p_{n}}-\sqrt{p_{n-1}}<B/2,\) where \( B\geq 4\) and \(\sqrt[4]{P_{A}}>B\), then between \(\left(x-1\right)^{2.5}\)\textit{ }and\textit{ }\(  x^{2.5}\),\textit{ }where \textit{x}\( \) is a real number greater than \( p_{A}\), there is always a prime number.

Theorem 10. Let there be a positive  integer constant \textit{A}, such that for any integer \( n\) greater than \( A\) consecutive primes \( p_{n-1},  p_{n}\) satisfy the inequality \(\sqrt{p_{n}}-\sqrt{p_{n-1}}<1,\) then between \(\left(x-1.5\sqrt[4]{x}\right)^{1.5}\)and\textit{ }\(  x^{1.5}\),\textit{ }where \textit{x}\( \) is a real number greater than \( p_{A}\), there is always a prime.

\section{The Legendre – Schinzel hypothesis and the Ribenboim concept}

According to W. Sierpiński: "\textit{It was Legendre who formulated the conjecture that for sufficiently large numbers x there is at least one prime between x} \textit{and}  \( x+\sqrt{x}\)" ([4], p.155). One of the possible versions of Legendre's statement was formulated by A. Schinzel ([4], p.155), which can be presented as follows:

Conjecture 1. There is a positive constant \textit{C}, such that for any real number \textit{x} greater than \textit{C} the open interval of real numbers \( (x,x+\sqrt{x})\) contains at least one prime number. 

Below we will show that the above hypothesis can be presented in a fundamentally different alternative form, in which only prime numbers is used.

Theorem 11. There is a positive constant \( C\), such that for any real number \textit{x} greater than \( C\) the open interval of real numbers \( (x,x+\sqrt{x})\) contains a prime if and only if there is a positive integer  constant\textit{ N}, such that for any integer \textit{n} greater than \textit{N} the open interval of real numbers \( (p_{n},p_{n}+\sqrt{p_{n}})\) contains a prime number.
\begin{proof}
Let Conjecture 1 be true then there is a positive integer \textit{N} such that \( p_{N-1}<C\leq p_{N}\) therefore, for any integer \textit{n} greater than \textit{N}, the open interval \(\left(p_{n},p_{n}+\sqrt{p_{n}}\right)\) contains at least one prime number \( q\). And vice versa also is true. Let there be a positive integer constant \textit{N} such that for any integer \textit{n} greater than \textit{N}, the open interval \(\left(p_{n},p_{n}+\sqrt{p_{n}}\right)\) contains at least one prime \( q\), then for any real number \textit{x} greater than \( p_{N+1}\), the open interval \(\left(x,x+\sqrt{x}\right)\) contains at least one prime number \textit{p}. Indeed, let there be some real number \( x_{0}\) greater than \( p_{N+1}\) such that the open interval \(\left(x_{0}, x_{0}+\sqrt{ x_{0}}\right)\) does not contain primes. Let us construct an open interval \(\left(p_{k},p_{k}+\sqrt{p_{k}}\right)\), where \( p_{N}<p_{k}<x_{0}<p_{k+1}\) since \( k>N\), therefore, \(\left(p_{k},p_{k}+\sqrt{p_{k}}\right)\) contains a prime. However,  on the other hand, the open interval \(\left(p_{k},p_{k}+\sqrt{p_{k}}\right) =\left(p_{k}, x_{0}\right]\cup\left(x_{0},p_{k}+\sqrt{p_{k}}\right)\), where \( x_{0}<p_{k}+\sqrt{p_{k}}\) otherwise \( p_{k}<p_{k+1}< x_{0}\), \( (x_{0},p_{k}+\sqrt{p_{k}})\subset (x_{0}, x_{0}+\sqrt{ x_{0}})\), since \( p_{k}+\sqrt{p_{k}}<x_{0}+\sqrt{ x_{0}}\), does not contain prime numbers. We have obtained a contradiction, and Theorem 11 is proved.
\end{proof}

Comment. Using the Legendre\textbf{ – }Schinzel hypothesis, which assumes that the difference of consecutive primes \( p_{n},  p_{n+1}\) satisfies the inequality \( p_{n+1}-p_{n}<\sqrt{p_{n}}\), it is not difficult to show that there exists a constant \( C\) such that for any integer \textit{m} greater than \( C\),  intervals \( (m^{2},\left(m+\frac{1}{2}\right)^{2})\) and \( (\left(m+\frac{1}{2}\right)^{2},\left(m+1\right)^{2})\) contain primes, we got exactly this result, guided by the Ribenboim concept and using the algorithmic proof model, Theorem 4, where we consider that primes \( p_{n},  p_{n+1}\) satisfy the inequality \(\sqrt{p_{n+1}}-\sqrt{p_{n}}<\frac{1}{2}\) and, as one of the consequences, the difference of primes \( p_{n},  p_{n+1}\) satisfies the inequality \( p_{n+1}-p_{n}<\sqrt{p_{n+1}}\).

\section{Prospects for using the Ribenboim concept in solving Problems with Prime Numbers}

Let's look at another possibility of using the Ribenboim concept: if it is true that there are integer positive constants \textit{A} and \textit{B} such that for any integer \( n\) greater than \textit{A}, any pair of consecutive prime numbers \( p_{n-1},  p_{n}\) satisfies the inequality \(\sqrt{p_{n}}-\sqrt{p_{n-1}}<\frac{1}{2B}\), which is consistent with a stronger version of estimating the difference of consecutive primes from the relation \(\sqrt{p_{n}}-\sqrt{p_{n-1}} =\) \( o(1)\), which is considered correct by researchers, for more details see [1], p.191, then the following statements are true:

Theorem 12. Let there be positive integer constants \textit{A} and \textit{B}, such that for any \textit{n }greater than \textit{A} the prime numbers \( p_{n-1},  p_{n}\) satisfy the inequality \(\sqrt{p_{n}}-\sqrt{p_{n-1}}<\frac{1}{2B}\), then for any integer \( m-1\) greater than \( p_{A}\), open intervals \( ((m-\frac{s+1}{2B})^{2},(m-\frac{s}{2B})^{2})\), where the integer \( s\) takes  values from \( 0\) to \( 2B-1\){\small  }, \( m>B\) and \(m \equiv0\mod{B}\), contain prime numbers.

Corollary 1. Let there be a positive integer constant \textit{A}, such that for any \textit{n }greater than \textit{A} the consecutive prime numbers \( p_{n-1},  p_{n}\) satisfy the inequality \(\sqrt{p_{n}}-\sqrt{p_{n-1}}<\frac{1}{2}\), then for any integer \( m-1\) greater than \( p_{A}\), the open intervals \( ((m-\frac{s+1}{2})^{2},(m-\frac{s}{2})^{2})\), where the integer \( s\) takes the values \( 0\); 1, contain prime numbers.

Theorem 13. Let there be positive integer constants \textit{A} and \textit{B}, such that for any \textit{n }greater than \textit{A} the prime numbers \( p_{n-1},  p_{n}\) satisfy the inequality \(\sqrt{p_{n}}-\sqrt{p_{n-1}}<\frac{1}{2B}\), then for any integer \( m-1\){\small  } greater than \( p_{A}\), open intervals \( ((m-\frac{s+1}{2B})^{2},(m-\frac{s}{2B})^{2})\), where the integer \( s\) takes values from \( 0\) to \( 2\sqrt{B}-1\){\small  }, \( m>B\) and \(m \not\equiv0\mod{B}\), contain prime numbers.

\section{Conclusion}

In this note, we proposed an original method of proving some statements concerning problems with prime numbers, so we proved the Ribenboim conjecture, in other words, we showed that this hypothesis really solves the Landau problem \# 4 (1912), from this proof it is clear that using the inequality \(\sqrt{p_{n+ 1}}-\sqrt{p_{n}} <1\) is crucial and cannot be weakened.

Next, we showed that if, starting from some integer \(N,\) consecutive prime numbers \(p_{n}\), \(p_{n+1}\) satisfy the inequality \(\sqrt{p_{n+1}}-\sqrt{p_{n}} <1/2\), then there exists an integer positive constant \textit{M} such that for any integer \(m\) greater than \textit{M} between the squares of any two numbers \(m\), \(m+\frac{1}{2}\) and \(m+\frac{1}{2}\), \( m +1\), there are always prime numbers.

Applying Ribenboim's approach, we obtained an original solution of Ingham's problem about primes between cubes of consecutive integers, which agrees with the results of [5].

We also proved that the solution of the Landau problem without the Ribenboim hypothesis, as a consequence, evaluates the difference of consecutive prime numbers \(p_{n}, p_{n+1}\) by the  relation \(p_{n+1}-p_{n} = O(\sqrt{p_{n}}\){\small)}.

We proved that the Legendre-Schinzel hypothesis (1961) about the existence of a prime number between the real numbers \(x\) and \(x+\sqrt{x}\) is equivalent to the hypothesis that the difference of consecutive prime numbers \(p_{n}, p_{n+1}\) satisfies the inequality \(p_{n+1}-p_{n}<\sqrt{p_{n}}\).

Thus, we have demonstrated a number of applications of the Ribenboim concept to known problems with prime numbers, both already solved and not yet solved, using  really expected estimates of the difference of consecutive prime numbers \( p_{n},  p_{n+1}\) from the inequality \(\sqrt{p_{n+1}}-\sqrt{p_{n}}<A\), where the values \( A\) are \(1/2;1;2\).

\begin{center}
\end{center}
\begin{center}
\(\) \(\) \(\) \(\) \(\) \(\) \(\) \(\) \(\)  REFERENCES
\end{center}

1.	P. Ribenboim, \(The\) \(Little\) \(Book\) \(of\) \( Bigger\) \( Primes\), Springer-Verlag, New York, 2004.

2.	J. Pintz, Landau's problems on primes, \(J.\) \(Th\acute{e}or.\) \(Nombres\) \(Bordeaux\), 21, 2, 357-404, 2009.

3.	 H. Iwaniec, M. Jutila, Primes in short intervals, \(Arkiv\) \(Mat\). 17,167–176, 1979.

4.	W. Sierpiński, \( Elementary\) \(Theory\) \(of\) \(Numbers,\) North-Holland, Amsterdam, 1988.

5.	A. Ingham, On the difference between consecutive primes, \(Quart\). \(J\). \(Math\). \(Oxford\), 8, 255-266, 1937.

\begin{center}
\end{center}
\begin{center}
\end{center}
Department of Mathematics and Statistics,\(\) \(\) \(\) \(\) \(\) \(\) \(\) \(\) \(\)\(\) \(\) \(\) \(\) \(\) \(\) \(\) \(\) \(\)\(\) \(\) \(\) \(\) \(\)\(\) \(\) \(\)\(\) \(\) \(\) \(\) \(\) \(\) \(\) \(\)\(\) \(\) \(\)\(\) \(\)\(\) \(\) \(\)\(\) \(\)\(\)\(\) \(\) \(\) \(\) \(\) \(\) \(\) \(\) \(\)\(\) \(\) \(\) \(\) \(\) \(\) \(\) \(\) \(\)\(\) \(\) \(\) \(\) \(\)\(\) \(\) \(\)\(\) \(\) \(\)  \(\) \(\)\(\) \(\) \(\) \(\) \(\)\(\) \(\) \(\)   \(\) \(\)\(\) \(\) \(\) 
Concordia University,\(\) \(\) \(\) \(\) \(\) \(\) \(\) \(\) \(\)\(\) \(\) \(\) \(\) \(\) \(\) \(\) \(\) \(\)\(\) \(\) \(\) \(\) \(\)\(\) \(\) \(\)\(\) \(\) \(\) \(\) \(\) \(\) \(\) \(\)\(\) \(\) \(\)\(\) \(\)\(\) \(\) \(\)\(\) \(\)\(\)\(\) \(\) \(\) \(\) \(\) \(\) \(\) \(\) \(\)\(\) \(\) \(\) \(\) \(\) \(\) \(\) \(\) \(\)\(\) \(\) \(\) \(\) \(\)\(\) \(\) \(\)\(\) \(\) \(\) \(\) \(\) \(\) \(\) \(\)\(\) \(\) \(\)\(\) \(\)\(\) \(\) \(\)\(\) \(\)\(\)  \(\)\(\) \(\) \(\)\(\) \(\)\(\)    \(\)\(\) \(\) \(\)\(\) \(\)\(\) \(\)\(\) \(\) \(\)\(\) \(\)\(\)  \(\)\(\) \(\) \(\)\(\) \(\)\(\) \(\)\(\) \(\) \(\)\(\) \(\)\(\) \(\) \(\)\(\) \(\) \(\) \(\) \(\)\(\) \(\) \(\)  \(\) \(\)\(\) \(\) \(\)         
1455 de Maisonneuve Blvd. W.,\(\) \(\) \(\) \(\) \(\) \(\) \(\) \(\) \(\)\(\) \(\) \(\) \(\) \(\) \(\) \(\) \(\) \(\)\(\) \(\) \(\) \(\) \(\)\(\) \(\) \(\)\(\) \(\) \(\) \(\) \(\) \(\) \(\) \(\)\(\) \(\) \(\)\(\) \(\)\(\) \(\) \(\)\(\) \(\)\(\)\(\) \(\) \(\) \(\) \(\) \(\) \(\) \(\) \(\)\(\) \(\) \(\) \(\) \(\) \(\) \(\) \(\) \(\)\(\) \(\) \(\) \(\) \(\)\(\) \(\) \(\)\(\) \(\) \(\) \(\) \(\) \(\) \(\) \(\)\(\) \(\) \(\)\(\) \(\)\(\) \(\) \(\)\(\) \(\)\(\) \(\) \(\) \(\) \(\) \(\) \(\)  
Montreal, QC H3G 1M8, Canada

\bibliography{references}
\bibliographystyle{acm}

\end{document}